\documentclass{article}
\usepackage[utf8]{inputenc}
\usepackage{latexsym,amsmath,amsthm,amssymb,graphicx,MnSymbol,calc}
\usepackage{tikz}\usetikzlibrary{matrix,arrows,calc,decorations,arrows.meta,positioning,shapes}
\usepackage[font=itshape]{quoting}
\usepackage{verse}
\usepackage{standalone}
\usepackage{float}
\theoremstyle{plain}
\newtheorem{thm}{Theorem}[section]

\theoremstyle{definition}
\newtheorem{defn}[thm]{Definition}

\newtheorem{exmp}[thm]{Example}

\newtheorem{rem}[thm]{Remark}
\usepackage[utf8]{inputenc}
\usepackage{graphicx}
\usepackage{enumitem}

\usepackage{latexsym,amsmath,amsthm,amssymb,graphicx,MnSymbol,calc}
\usepackage{tikz}\usetikzlibrary{matrix,arrows,calc,decorations,arrows.meta,positioning,shapes}
\usepackage{graphicx}
\definecolor{woured}{RGB}{227,24,55}
\definecolor{woublue}{RGB}{24,85,227}

\title{Multigraphs from crossword puzzle grid designs}
\author{Ben Cot\'e \\ Western Oregon University \and Leanne Merrill \\ Western Oregon University}
\usepackage{graphicx}
\begin{document}

\maketitle

\begin{abstract}
    
Crossword puzzles lend themselves to mathematical inquiry. Several authors have already described the arrangement of crossword grids and associated combinatorics of answer numbers \cite{ferland2014record} \cite{ferland2020counting} \cite{McSweeny2016Analysis}. In this paper, we present a new graph-theoretic representation of crossword puzzle grid designs and describe the mathematical conditions placed on these graphs by well-known crossword construction conventions.
\end{abstract}

Solving a crossword puzzle is not an inherently mathematical task, but their underlying structure raises many interesting mathematical questions. For instance: how many allowable crossword puzzle grids of a given size are there? This paper provides new tools to approach such questions.

In Section ~\ref{section1}, we give the basic definitions and rules for the class of crossword puzzles under consideration. In Section~\ref{section2}, we describe a graph associated with a given crossword puzzle grid. In Section~\ref{section3}, we define a class of graphs that includes those graphs associated to puzzles. Section~\ref{section4} contains a set of necessary conditions for these graphs to represent crossword puzzles. Finally, in Section~\ref{section5}, we describe the next phase of the project that will allow us to count crossword puzzle grids using this special class of graphs via combinatorial arguments.

\section{Background}\label{section1}

We start with a definition adapted from \cite{McSweeny2016Analysis}: 

\begin{defn}[Crossword puzzle]
A \emph{crossword puzzle} is a grid made up of squares that are either white or black.  The white squares are called \emph{cells} and the black squares are called \emph{voids}. A maximal vertical sequence of adjacent cells is a \emph{Down answer}, and likewise a maximal horizontal sequence of cells is called an \emph{Across answer}. The \emph{boundary} of a puzzle refers to perimeter of the crossword grid.\end{defn}

Many types of crossword puzzles exist, but the most commonly published in US newspapers are \emph{American-style} crossword puzzles.  These puzzles consist of a  $(2n+1) \times (2n+1)$ grid for a natural number $n$ obeying certain \emph{structure rules} (see \cite{BasicRul24:online}):  

\begin{enumerate}
    \item (Connectivity) Any two cells in the puzzle are connected via a path traveling only through cells using horizontal or vertical paths. In other words, the voids in a puzzle do not ``break up" the puzzle into smaller, disconnected portions; 
    \item ($180^\circ$ rotational symmetry) The grid must possess rotational symmetry. That is, when it is rotated $180^\circ$, the grid must look the same; 
    \item (Answer length) All answers must be at least three cells long; 
    \item (Keyed squares) Each cell must \emph{keyed}; that is, each cell must be part of both a vertical and horizontal answer;
    \item (Full dimension) There cannot be a full row or column of voids along the boundary of the puzzle. 
\end{enumerate}

Both McSweeney \cite{McSweeny2016Analysis} and Ferland \cite{ferland2014record}, \cite{ferland2020counting} describe crossword puzzles mathematically. In \cite{McSweeny2016Analysis}, a bipartite graph associated with a more generalized crossword construction is given. In \cite{ferland2014record}, Ferland counts the maximum number of clues in a standard American-style crossword puzzle by exploiting the structure rules, particularly symmetry and answer length. In the next section, we combine and modify these approaches to define a new type of graph that completely captures the structure of a crossword puzzle grid.

\section{The crossword multigraph}\label{section2}

To each crossword puzzle, there is an associated \emph{crossword network} with a vertex for each answer and an edge between the vertices associated to any two answers that intersect in a cell  \cite{McSweeny2016Analysis}.  Since the set of answers is partitioned into Down answers (the set $Down$) and Across answers (the set $Across$), and no two answers of the same type can intersect, the crossword network is a bipartite graph.  Additionally, the edges correspond uniquely to the cells of the puzzle, so degrees of vertices give the length of each answer. Connectivity in the puzzle is precisely connectivity in the graph. 

From a crossword network, however, it would be difficult to reconstruct the crossword grid to which it belongs. Since the network is bipartite, one could easily partition the vertices into two sets, but it is usually impossible to determine which set is Down, which is Across, or the relative arrangement of vertices (answers) in each part.   

On the other hand, the crossword network contains duplicate information. Due to the $180^\circ$ rotational symmetry, answers which do not contain the center cell come in pairs. For the crossword network, this pairing corresponds to an order two graph automorphism which respects the partition of vertices into $Down$ and $Across$. If we are only interested in the grid design, there is no need to distinguish between symmetric pairs of answers. 

Our first goal is to arrive at a graph which addresses these issues.  To start, we introduce coordinates for each cell.

\begin{figure}[ht]
    \centering
    \includegraphics[]{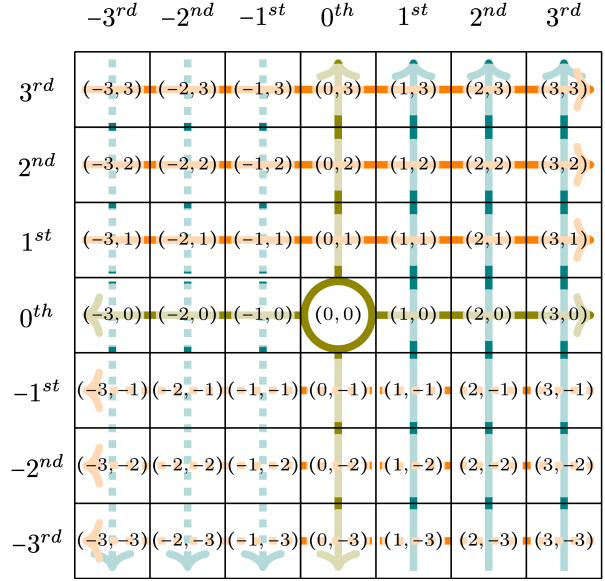}
    \caption{Cell, column, and row labeling with arrows denoting row and column orientation}
    \label{fig:cell_label}
\end{figure}

\begin{defn}[Puzzle coordinates and row/column orientation]
 For a $(2n+1)\times (2n+1)$ puzzle, the \emph{coordinates} of each cell are determined by assigning the center cell the pair $(0,0)$ and expanding out with typical Cartesian conventions.  With this, the \emph{$i^{th}$ column} will be the set of cells of the form $(i,k)$ with $-n \leq k \leq n$ and rows are analogously defined. The cell in the $i^{th}$ column and $j^{th}$ row is labeled $(i,j)$. Further, each row and column are given an \emph{orientation}. Zero rows and columns are oriented outwards from the origin. Positive rows are oriented from left to right and positive columns are oriented from bottom to top.  Negative rows and columns are given an orientation opposite that of their positive counterparts. See Figure~\ref{fig:cell_label} for an example, and notice that the orientations are preserved by the rotational symmetry.
\end{defn}

 Using the row and column numbers and the orientation, we can assign numerical indices to the vertices and signed labels to the edges of the crossword network.

\begin{figure}
    \centering
\includegraphics[]{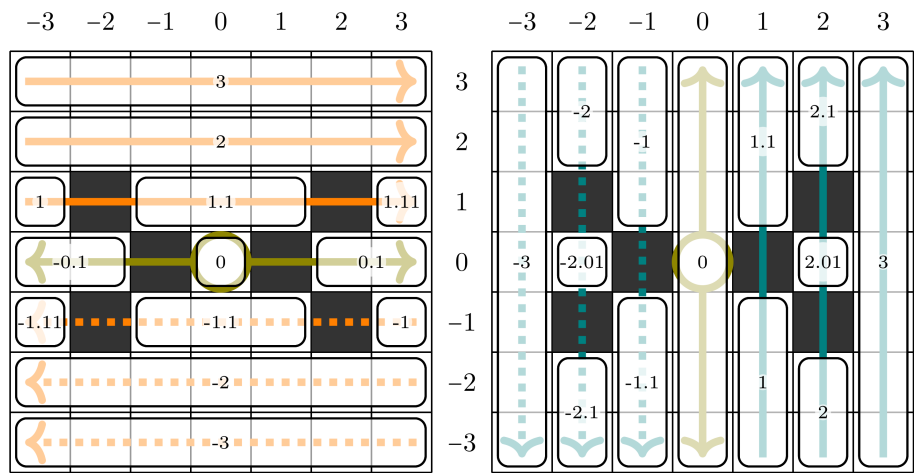}
            \caption{A puzzle showing an indexing of Across (left) and Down (right).  Answers are indexed so their absolute values increase with the orientation of the row/column and the indices negate under the rotational symmetry.}
    \label{fig:orientation}
\end{figure}

\begin{figure}
    \centering
        \includegraphics[]{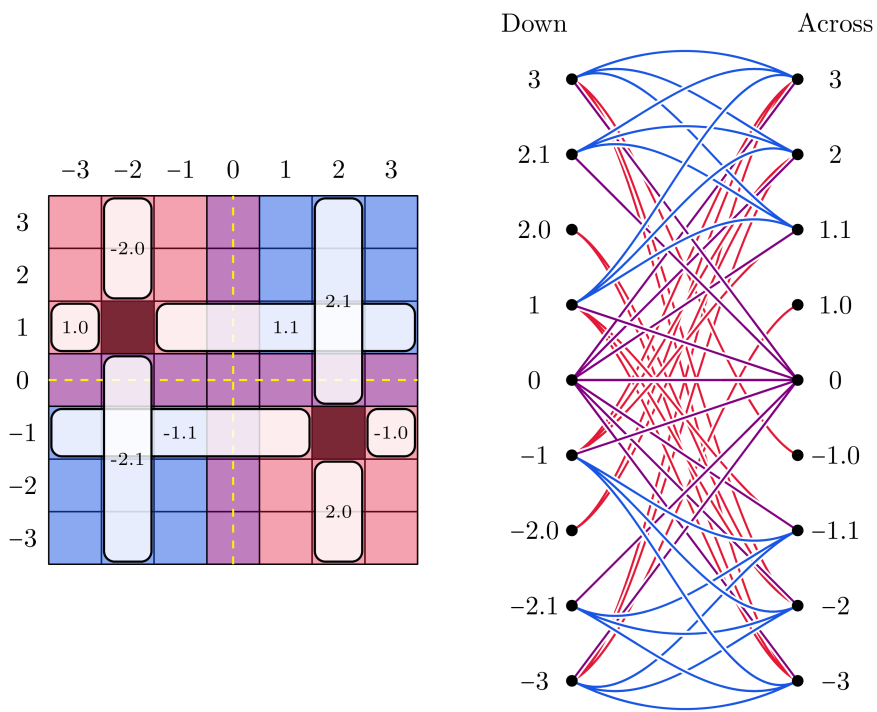}
    \caption{A crossword puzzle and the corresponding labeled indexed crossword network. Note that this puzzle does not obey the answer length rule.}
    \label{fig:indexednetwork}
\end{figure}

\begin{defn}[Labeled indexed crossword network or LICN]
   Each answer is \emph{indexed} with a rational decimal number where the integer part of the index represents the row or column number.  When a row or column is comprised of multiple answers, the indices also contain decimal parts, the absolute values of which increase along the direction of the orientation.  Further, we require these indices to be negated under the $180^\circ$ rotational symmetry.  See Figure~\ref{fig:orientation} for an example of such an indexing.  Each edge is assigned a \emph{label} based on the product of its coordinates: ``$-$" if negative, ``$0$" if $0$, and ``$+$" if positive.  

   The crossword network, when indexed and labeled in this way, we call the \emph{labeled indexed crossword network}. To avoid clutter in our graphs, we use colors to indicate the labels of the edges, with ``$+$" colored blue, ``$-$" colored red, and ``$0$" colored purple. See Figure~\ref{fig:indexednetwork}.

\end{defn}

\begin{figure}[ht]
    \centering
        \includegraphics[]{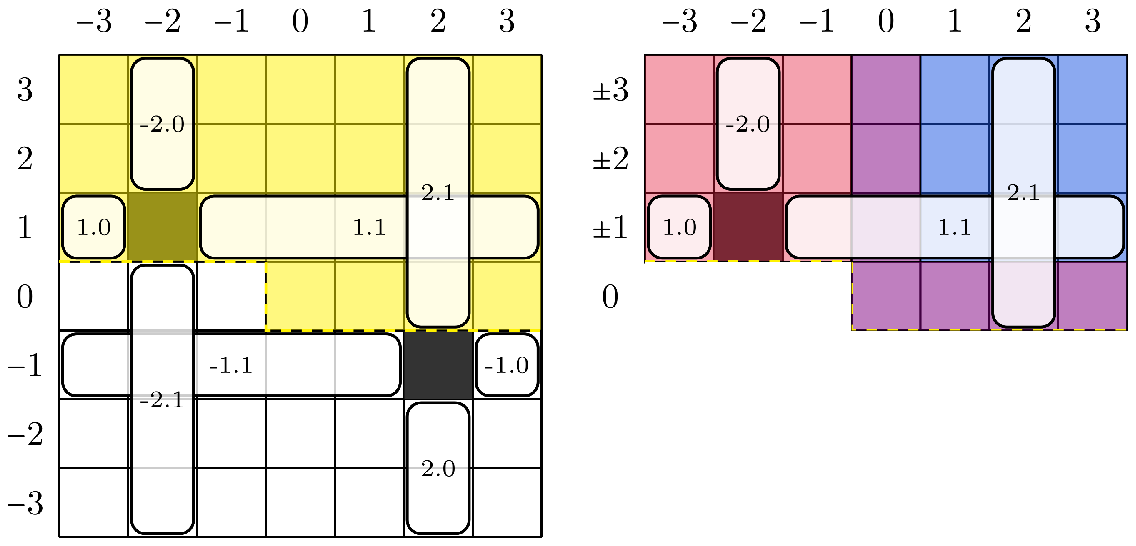}
    \caption{On the left, the fundamental region shaded in yellow.  On the right, the fold is shown by the red and blue halves of the fundamental region.}
    \label{fig:fundreg}
\end{figure}

Notice that the symmetry in the crossword puzzle leads to a symmetry in the labeled indexed crossword network.  In particular, the vertices with positive index have a corresponding vertex with negative index.  In order to reconstruct the grid, we would only need one of each of these pairs. We choose the positive rows, along with the nonnegative part of the $0^{th}$ row as a fundamental region for the $180^\circ$ rotational symmetry (see Figure~\ref{fig:fundreg}).

For ease of demonstration, we provide a $5 \times 5$ crossword grid design and the corresponding fundamental region in Figure~\ref{fig:fundreg}, though it breaks the answer length rule. We will continue to use this example throughout the section.

To the fundamental region, we associate a graph derived from the labeled indexed crossword network. 

\begin{defn}[Fundamental graph]
The \emph{fundamental graph} is subgraph of the LICN when restricted to the fundamental region of the grid.  The indices of the nonzero vertices in Across are marked with a $\pm$ since they now designate a positive and negative index from the original LICN.  
\end{defn}

The choice of taking half the rows or half the columns for a fundamental region leads to two ways to fold the graph, one that collapses pairs of symmetric Across answers and one that collapses pairs of symmetric Down answers; we choose to collapse the Across answers.

\begin{figure}[ht!]
    \centering
        \includegraphics[]{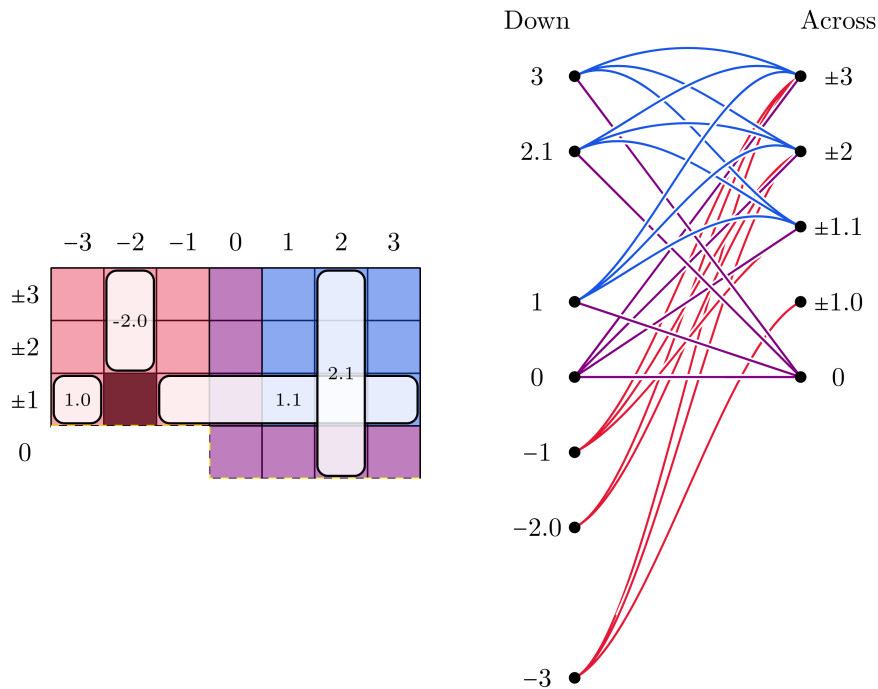}
    \caption{On the left, the fundamental region.  On the right, the corresponding fundamental graph.}
    \label{fig:fungraph}
\end{figure}

The fundamental graph of the crossword puzzle is not the ultimate goal of our study, but it is an important step that uses the symmetry of the puzzle to simplify our graphical representation. In the fundamental graph, it is no longer the case that the degree of a vertex corresponds to the length of an answer. In our example, since the answer $3$ Down is not contained within the fundamental region, its cells are now represented by both the edges from $3$ Down and from $-3$ Down.

We remedy this issue, and create a much more compact representation, with our main definition for this section.

\begin{figure}[h]
    \centering
        \includegraphics[]{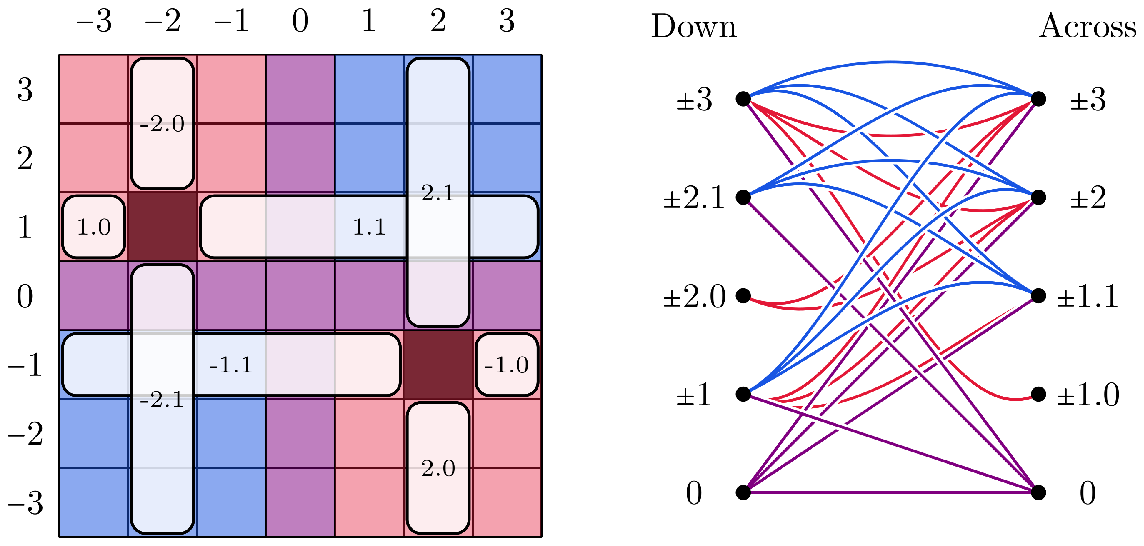}
    \caption{On the left, a crossword grid with only two voids.  On the right, the corresponding crossword multigraph.}
    \label{fig:crossmult}
\end{figure}

\begin{defn}[Crossword multigraph]
Given a $(2n+1)\times (2n+1)$ crossword grid, the \emph{crossword multigraph} is the multigraph formed from the fundamental graph by identifying vertices from Down that have equal absolute values.  For example, the Down vertices $3$ and $-3$ in Figure~\ref{fig:fungraph} become the Down vertex $\pm3$ in Figure~\ref{fig:crossmult}. The result is a multigraph where a pair of vertices might be adjacent through both a blue and red edge (e.g. $\pm 3$ Down and $\pm2$ Across in Figure~\ref{fig:crossmult}).
\end{defn}

In the crossword multigraph, we regain the equality between degree and word length, though each nonzero vertex now represents two (symmetric) answers instead of the single answer represented in the original crossword network. We also see that nearly every edge corresponds to two cells in the original puzzle (one from the fundamental domain and one outside of it), though the single $0-0$ edge, if it exists, corresponds uniquely to the center cell of the puzzle if it is indeed a cell and not a void.

From a crossword multigraph one can reconstruct the crossword grid.  In particular, for any edge one can determine the cell(s) to which it corresponds by taking the integer parts of the indices of its incident vertices as coordinates, and applying $+$ or $-$ signs so the coordinates multiply to the edge's label.  Additionally, the crossword multigraph is essentially determined by the crossword grid.  The only choices one can make along the way are the decimal parts to use in the indexing.  Together, these observations prove the following:

\begin{thm}[One-to-one correspondence]
There is a unique crossword multigraph for each crossword puzzle design, up to order-preserving reindexing of the vertices.
\end{thm}

We wish to investigate when an arbitrary graph is a crossword multigraph, and if so, its corresponding crossword puzzle grid shape.  In order to do so, we first need to define a general class of graphs which contain the basic structure, indices, and labeling.

\section{Bit multigraphs and voiding} \label{section3}

We now introduce a class of graphs that we will study in detail. 

\begin{defn}[Bipartite indexed tricolored (bit) multigraph]

Let $(V,E)$ be a graph. We say that $(V,E)$ is a \emph{bit multigraph} if it satisfies the following conditions: 

\begin{itemize}
    \item \textbf{Balanced Bipartite:} the graph is bipartite, consisting of two vertex sets $A$ and $B$ such that $A \cup B = V$ and $|A| = |B|$;
    
    \item \textbf{Indexed:} the vertices within each part are indexed with distinct nonnegative decimal numbers.  Vertices whose index has zero integer part are called \emph{zero vertices} and those with non-zero integer part are called \emph{non-zero vertices}. Additionally, nonzero vertices are decorated with a $\pm$ before the decimal index. $A$ and $B$ must each contain at least one zero vertex (which may be isolated).  
    
    \item \textbf{Tricolored} each edge is given a label from the set $\{-,0,+\}$ satisfying the following restrictions:
        \begin{enumerate}
            \item If a pair of adjacent vertices includes a a zero vertex, then they must be adjacent with only one edge and it must be labeled $0$.  Only one pair of vertices with zero integer part can be adjacent.
            \item All other edges are labeled $-$ or $+$, and a pair of vertices may be adjacent through both a $-$ edge and a $+$ edge.
        \end{enumerate}
     Visually, we distinguish between edge labeling by coloring the edges blue if labeled $+$, purple if labeled $0$, and red if labeled $-$.  

\end{itemize}

\end{defn}

A consequence of the restrictions in the definition above is that two vertices may have at most two edges between them. If indeed a pair of vertices are adjacent by two edges, those edges must have labels $+$ and $-$. 

We have seen such objects before. 
\begin{figure}
    \centering
    \includegraphics[]{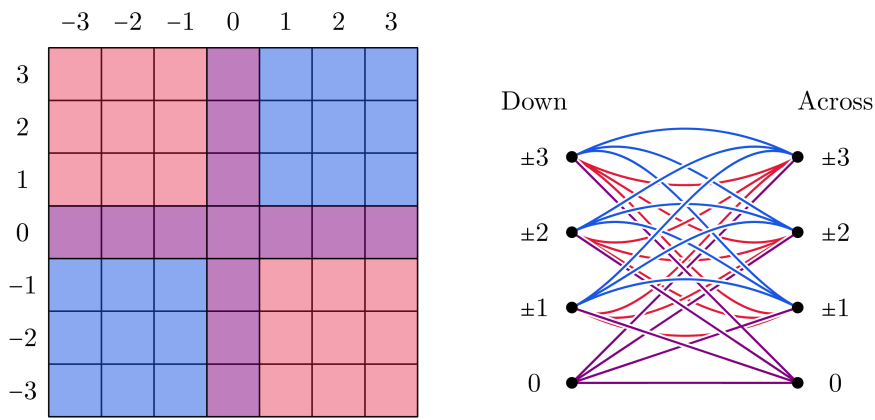}
    \caption{Crossword with no voids and corresponding bit multigraph.}
    \label{fig:novoid}
\end{figure}

\begin{exmp}[Crossword multigraph with no voids]\label{exmp:unvoidedbit}
The crossword multigraph of a $(2n+1)\times(2n+1)$ crossword grid with no voids is a bit multigraph with both parts containing vertices indexed $\{0,\pm1,\pm2,\ldots \pm n\}$.  In particular, it is a $K_{n,n}$ of edges labeled $+$, a $K_{n,n}$ of edges labeled $-$, and a pair of $K_{1,n}$ graphs with their degree $n$ vertices connected, with all of these $2n+1$ edges labeled $0$. All other pairs of adjacent vertices are adjacent by both a $+$ and $-$ labeled edge. See Figure~\ref{fig:novoid} for a $7 \times 7$ example. 
\end{exmp}

We now describe a \emph{voiding} procedure which takes in a given bit multigraph and produces another bit multigraph with two additional vertices and one fewer edge than the original graph. 

\begin{defn}[Voiding procedure]
Given a bit multigraph, an edge can be \emph{voided} by performing the following steps:
\begin{enumerate}
    \item The edge is removed. 
    \item Each vertex incident to removed edge is ``split" by removing that vertex and replacing that vertex with two new vertices.  The indices of these new vertices are the same as the vertex they came from but with an additional trailing decimal place.  One is given the decimal $0$ and the other the decimal $1$ in that place. We will use the notation that a vertex indexed $a$ splits into $a0$ and $a1$.  
    \item \label{itm:void_assign_edges} Suppose the indices of the removed vertices were $a$ and $b$ respectively.  The remaining edges incident to each split vertex are assigned to the new vertices using the following rules: 
    \begin{enumerate}
        \item \label{itm:void_labels_match} Edges whose label matches the label of the removed edge are assigned to the new vertices using the following rules:
    \begin{enumerate}
        \item If an edge $\{a,y\}$ is labeled $0$ or $+$, then it is assigned to vertex $a0$ if $y< b$, and assigned to $a1$ if $y> b$.  
        \item If an edge $\{a,y\}$ is labeled $-$, then it is assigned to vertex $a0$ if $y > b$ and assigned to $a1$ if $y < b$.  
         \item If an edge $\{x,b\}$ is labeled $0$ or $+$, then it is assigned to vertex $b0$ if $x< a$, and assigned to $b1$ if $x> a$.
        \item If an edge $\{x,b\}$ is labeled $-$, then it is assigned to vertex $b0$ if $x > a$ and assigned to $b1$ if $x <  a$.  
        \end{enumerate}
        
    \item\label{itm:void_labels_dont_match} Edges whose labels do not match the label of the removed edge are assigned to the new vertices using the following rules: 
    \begin{enumerate}
        \item If an edge $\{a,y\}$ is labeled $+$, then it is assigned to vertex $a1$.  
        \item If an edge $\{a,y\}$ is labeled $0$, there are two cases: it is assigned to vertex $a1$ if the removed edge was labeled $-$, and assigned to vertex $a0$ if the removed edge was labeled $+$. 
        \item If an edge $\{a,y\}$ is labeled $-$, then it is assigned to vertex $a0$. 
        \item If an edge $\{x,b\}$ is labeled $+$, then it is assigned to vertex $b1$.  
        \item If an edge $\{x,b\}$ is labeled $0$, there are two cases: it is assigned to vertex $b1$ if the removed edge was labeled $-$, and assigned to $b0$ if the removed edge was labeled $+$. 
        \item If an edge $\{x,b\}$ is labeled $-$, then it is assigned to vertex $b0$.  
        \end{enumerate}
\end{enumerate}
\end{enumerate}

\end{defn}

\begin{figure}
    \centering
        \includegraphics[]{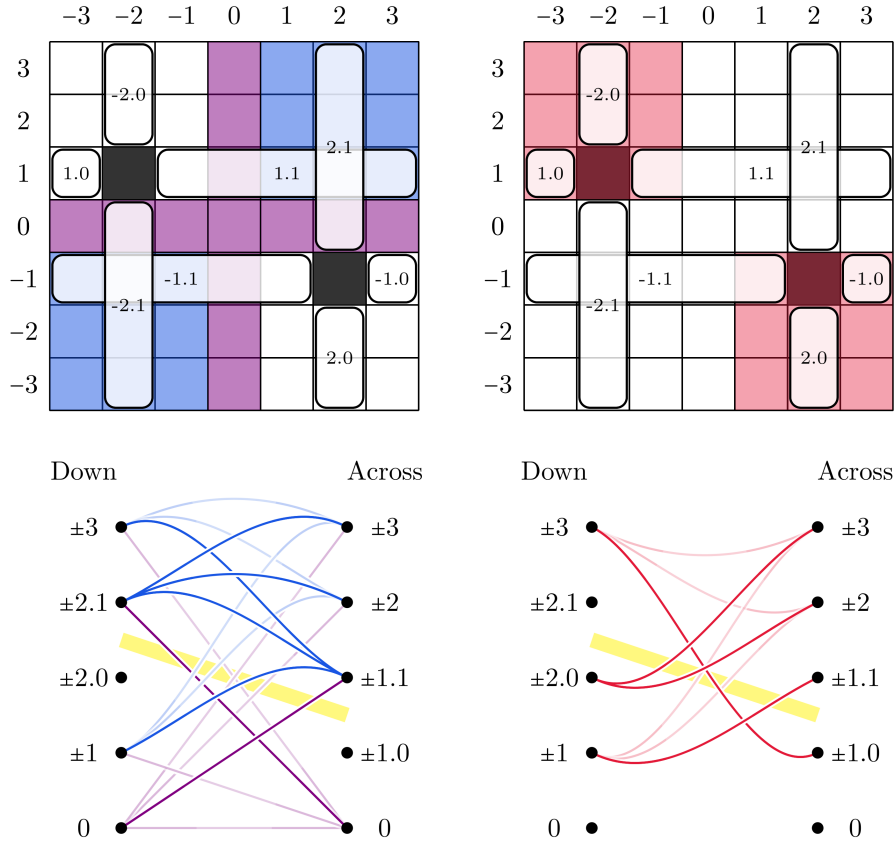}
    \caption{Example of voiding procedure emphasizing affected edges labeled $0$ or $+$ on the left, and $-$ on the right.}
    \label{fig:voidex}
\end{figure}

\begin{rem}\label{updown} 
It is helpful to visualize Step~(\ref{itm:void_assign_edges}) of the voiding procedure by describing how the affected edges distribute between the split vertices. First, note that in a bit multigraph, the indexing of the vertices allows the vertices to be totally ordered in a vertical direction from least to greatest in each part, and we adopt this convention with all bit multigraphs. The conditions in Step~(\ref{itm:void_labels_match}) describe when certain edges will cross the yellow ``void line" seen in Figure~\ref{fig:voidex}. Using the vertical ordering,  Step~(\ref{itm:void_labels_dont_match}) can be visualized as certain edges going ``up" to the vertex with index appended by a $1$ or ``down" to the vertex with index appended by a $0$, where ``up" and ``down" refer to the ordering of the indices in each part. 
\end{rem}

Using these ideas, the voiding procedure can be summarized with the following limerick:

\settowidth{\versewidth}{When the void and edge have the same signs:}
\begin{verse}[\versewidth]
When the void and edge have the same signs:\\
only red edges cross the void lines.\\
\vin Else: down with red, up with blue,\\
\vin purple depends on void's hue\\
(down if blue but if red it inclines).
\end{verse}

Since the index of each ``split" vertex is formed by appending binary digits to the index of the original vertex, the new indices introduced in this way will maintain their order relative to the existing indices. This means that the edge reassignments in the third step of the voiding procedure will be consistent no matter the order in which the voided edges are removed. 

As seen in  Example~\ref{exmp:unvoidedbit}, the crossword multigraph of an unvoided crossword grid is a bit multigraph.  Even more, the voiding procedure encodes the changes made to a crossword multigraph when a cell becomes a void. 

\begin{thm}[Crossword voids]\label{thm:crossvoid}
Given a crossword grid, starting with the unvoided bit multigraph of the same size (e.g. Figure~\ref{fig:novoid}) and applying the voiding procedure to the edge representative of each void in the fundamental region results in a bit multigraph isomorphic to the crossword multigraph with additional isolated vertices from every place where either two voids share a side or a void shares a side with the boundary.
\end{thm}

\begin{proof}
We proceed by induction on the number of voids within the fundamental region.  If there are no cells voided from the fundamental region, then as in Example~\ref{exmp:unvoidedbit} the graphs are isomorphic.  Suppose that the result holds for any crossword grid that has $k$ voids in the fundamental region.  We will analyze a crossword grid with $k+1$ voids in the fundamental region by considering the additional void.

If the void is neither adjacent to another void nor to a boundary, then the vertices incident to the edge removed are neither maximally indexed in their parts nor previously split.  The voiding procedure partitions the edges to exactly match the partition of the Across answer and Down answer to which the cell belonged.

If the void is adjacent to a boundary, then one of the vertices is the maximal non-isolated vertex in its part.  Thus, the voiding procedure will produce an isolated vertex above it.

If the void is adjacent to another void, then the edge corresponding to the cell to be voided is incident to a vertex which resulted from a previous voiding operation. When this vertex is split, its edges will all go to the vertex with last two decimal digits agreeing.  The other will be an isolated vertex. 

\end{proof}

\section{Properties of crossword bit multigraphs} \label{section4}

In this section, we establish a set of necessary conditions for a bit multigraph to correspond to an allowable crossword puzzle grid (that is, one following the structure rules in Section~\ref{section1}). We introduce a definition that helps us keep track of the vertex indices in a given bit multigraph. 

\newcommand{\afloor}[1]{\lfloor #1 \rfloor_A}
\newcommand{\bfloor}[1]{\lfloor #1 \rfloor_B}
\newcommand{\floor}[1]{\lfloor #1 \rfloor}
\begin{defn}[Floor sets]
 The vertex indices partition $Across$ and $Down$ into \emph{floor sets} by the floor function applied to the indices.  For example, $\afloor{3}$ denotes the set of vertices in the part $A$ that are indexed with a number $3\leq j<4$.  
\end{defn}

As a notation reminder, recall that for a set of vertices $V$, let $E(V)$ denote the set of all edges incident to a vertex in $V$. Additionally, recall Remark~\ref{updown} for the vocabulary relating to edges going ``up" or ``down." 

Now we are ready to state and prove the main theorem of this paper. 

\begin{thm}[Necessary conditions on bit multigraphs to represent crossword puzzles.]\label{thm:listprop} Let $(V,E)$ be a bit multigraph with parts $A$ and $B$ as its vertex sets. If $(V,E)$ is the bit multigraph corresponding to a $(2n+1) \times (2n+1)$ crossword puzzle $C$, then: 
\begin{enumerate}
    \item (Squareness) Both $A$ and $B$ are partitioned into $n+1$ floor sets.\label{floorsetpart}  
    \item (Word length) Every non-isolated vertex has degree $\geq3$ \label{wordlength} 
    \item (Connectivity) The non-isolated subgraph is connected.\label{connectivity} 
    \item (Edge and vertex combinatorics) With $e=|E|$ edges and $k=|A|=|B|$ in each part, $e+k=2n^2+3n+2$.  That is, $e+k$ must be one more than a triangular number with an odd base (see \cite{oeis}).
    \label{triangle} 
    \item (Floor set edge conditions)  Here we discuss conditions on floor sets from $A$.  The conditions hold symmetrically for floor sets from $B$. 
    \begin{enumerate}
           \item (Zero floor set edge count) \label{thm:zero_floor} Let $|\afloor{0}|= l$, the number of vertices in $\afloor{0}$.  Then: %n+1-(l-1)  
          
           \begin{subequations}
           \begin{align}
            |E(\afloor{0})| &= n -l+2 \label{eq:zero_count} \tag{5.a.1} \\
                            &\geq |\{r : |\bfloor{r}|=1 \}| \label{eq:zero_count_ineq}   \tag{5.a.2}          
           \end{align}
        \end{subequations}
           
          \item (Nonzero floor set edge count) \label{thm:nonzero_floor} For $k>0$, $E(\afloor{k})$ contains at most $n$ edges labeled $-$, $n$ edges labeled $+$, and one edge labeled $0$. Setting $|\afloor{k}|= l$,  
           \begin{align}       
           |E(\afloor{k})| &= 2n -l+2 \label{eq:nonzero_count} \tag{5.b.1}\\
                        &\geq 2|\{r : |\bfloor{r}|=1 \}| \label{eq:nonzero_count_ineq} \tag{5.b.2}   
           \end{align} %2n+1-(m-1) 

         \item (No doubles) \label{thm:no_dubs} Between two nonzero floor sets, there can be at most one $+$ and one $-$ edge.  A zero floor set can be connected to any floor set by at most one $0$ edge.

        \item (Blue above red) \label{thm:void1} If $a_1$ and $a_2$ are in $\afloor{k}$ with $0<a_1<a_2$ and both $a_1$ and $a_2$ are incident to $b\in B$ with $b>0$, then edge $\{a_1,b\}\in -$ and edge $\{a_2,b\}\in +$. %  p.56 goodnotes  

        \item (Purple in between) \label{purp_inbetween}  If any vertex is incident to both an edge labeled $+$ and an edge labeled $-$, then it must be incident to an edge labeled $0$.

         \item (Maximal set of same labeled edges) \label{thm:max_edge} If a floor set has $n$ edges labeled $+$ or $n$ edges labeled $-$, then the set of same-labeled edges all must be incident to the same vertex in that floor set.  Further, if they are labeled $+$ they must be incident to the maximally indexed vertex and if labeled $-$ they must be incident to the minimally indexed vertex. 
         
         \item (Blue sweep) \label{thm:blue_sweep} Suppose there is a floor set $\afloor{k}$ which is adjacent via a $+$ edge to all of the floor sets in $B$ other than  $\bfloor{m_1},\bfloor{m_2},\ldots,\bfloor{m_l}$ with $m_1 > m_2 > \dots > m_l$. Then the $i^{th}$ highest index vertex in $\afloor{k}$ is adjacent to precisely the floor sets $\bfloor{x}$ with $m_{i-1}>x>m_{i}$, where $m_0=\infty$. Informally, this means that blue edges cannot cross and the $k^{th}$ highest floor set in $A$ can only be connected with edges labeled $+$ to the $k^{th}$ highest, or higher vertex within a given floor set in $B$. 
         
         \item (Red sweep) \label{thm:red_sweep} 
         Suppose there is a floor set $\afloor{k}$ which is adjacent via a $-$ edge to all of the floor sets in $B$ other than  $\bfloor{m_1},\bfloor{m_2},\ldots,\bfloor{m_l}$ with $m_1 > m_2 > \dots > m_l$. Then the $i^{th}$ lowest index vertex in $\afloor{k}$ is adjacent to precisely the floor sets $\bfloor{x}$ with $m_{i-1}>x>m_{i}$, where $m_0=\infty$. As a consequence, the $k^{th}$ highest floor set in $A$ can only be connected with edges labeled $-$ to the $k^{th}$ lowest, or lower vertex within a given floor set in $B$. 
         
         \item (Purple sweep) \label{thm:purp_sweep}
         Suppose the floor set $\afloor{0}$ is adjacent via a $0$ edge to all of the floor sets in $B$ other than  $\bfloor{m_1},\bfloor{m_2},\ldots,\bfloor{m_l}$ with $m_1 > m_2 > \dots > m_l$. Then the $i^{th}$ highest index vertex in $\afloor{0}$ is adjacent to precisely the floor sets $\bfloor{x}$ with $m_{i-1}>x>m_{i}$, where $m_0=\infty$. If $\afloor{0}$ and $\bfloor{0}$ are adjacent, then the connected vertices must be minimally indexed within each floor set.  If $\afloor{0}$ is adjacent to $x$, then no vertex in $\bfloor{x}$ of higher index than $x$ can be incident to a $-$ edge, and no vertex in $\bfloor{x}$ of lower index than $x$ can be incident to a $+$ edge. 
         
    \end{enumerate}
\end{enumerate}
\end{thm}

\begin{proof} 
The first three properties follow directly from the basic construction and structure rules for crossword puzzles. Property~\ref{floorsetpart} follows from the fact that each floor set in each part corresponds to a row or column in the fundamental region of the puzzle, and there are $n+1$ such rows and columns. 
The structure rules specify that each word length must be at least 3; word length in the crossword multigraph corresponds exactly to vertex degree, and thus Property~\ref{wordlength} must hold. Finally, Property~\ref{connectivity} follows exactly from the connectivity structure rule, allowing for isolated vertices resulting from the voiding procedure. 

Property~\ref{triangle} is the first property that does not follow directly from the structure rules. Let $v = k - (n+1)$. Since $k$ is the number of vertices in each floor part, and a crossword puzzle with no voids starts with $n+1$ floor sets, $v$ counts the number of voiding operations performed on a crossword multigraph with no voids (that is, the graph described in Example~\ref{exmp:unvoidedbit}) to obtain the puzzle $C$ represented by the given bit multigraph. To find $e$, note that in a crossword multigraph with no voids, there are $2n^2 + 2n + 1$ edges:  $n^2$ for each copy of $K_{n,n}$, and $2n +1$ for the edges labeled $0$. Each void in the fundamental region removes an edge; thus in a general bit multigraph, there are $e = 2n^2 + 2n +1 - v$ edges. 

Combining these equations, we get
\begin{eqnarray*}
e+k  &= &2n^2 + 2n + 1 - v + (n+1) + v \\ 
&=& 2n^2 + 3n + 2
\end{eqnarray*}
Such numbers are precisely one more than triangular numbers with an odd base, since 

$$2n^2 + 3n + 2 = \frac{(2n+1) ((2n+1) + 1)}{2} + 1$$

To prove \ref{thm:zero_floor}, consider the set $\afloor{0}$.  If the graph is unvoided, then $|\afloor{0}|= 1$ and $|E(\afloor{0})|=n+1$.  Whenever an edge incident to $\afloor{0}$ is voided, $|\afloor{0}|$ increases by $1$ and $|E(\afloor{0})|$ decreases by $1$.  After $l$ removals of incident edges, $|E(\afloor{0})|=n+1-(l-1)$. Additionally, every floor set in $B$ of size $1$ could not have been split by a voiding operation, so each must be adjacent to $\afloor{0}$. The inequality holds since each size 1 floor set needs to be adjacent to a zero vertex. For \ref{thm:nonzero_floor}, since an unvoided graph has precisely $n$ edges labeled $+$ and $n$ edges labeled $-$ emanating from a given nonzero vertex, and voiding removed edges, both the equality and inequality hold for nonzero vertices also. 

 For an unvoided crossword grid, \ref{thm:no_dubs} holds and is seen to be preserved by the voiding procedure since edges are only ever removed. By \ref{thm:no_dubs}, ${a_1,b}$ and ${a_2,b}$ cannot have the same label and cannot be labeled $0$.  By Theorem~\ref{thm:crossvoid}, the graph must be isomorphic to the multigraph obtained by the voiding procedure.  No edge between $\afloor{a}$ and $\bfloor{b}$ could have been voided, since both $+$ and $-$ edge are still present. The only way the edges from $\afloor{a}$ to $\bfloor{b}$ would split up is if a $+$ edge from $\afloor{a}$ was voided, in which case the $+$ edge goes ``up" and the $-$ edge goes ``down".  Thus \ref{thm:void1} is proven. 

Parts \ref{thm:max_edge}, \ref{purp_inbetween}, \ref{thm:blue_sweep}, \ref{thm:red_sweep}, and \ref{thm:purp_sweep} follow from a careful consideration of the voiding procedure, and in particular the different behavior of edges labeled $+$, edges labeled $-$, and edges labeled $0$  when a void is introduced.

\end{proof}

\section{Future directions}\label{section5}

We conjecture that the necessary conditions listed in Theorem~\ref{thm:listprop} are sufficient as well.   Once a bijection is finalized between crossword grid designs and bit multigraphs with particular properties, we can use the multigraphs in order to answer questions about the crossword grid designs. One such question is how many crossword grids exist of a given size.  For puzzles up to $15 \times 15$, these numbers were computed by Jim Ferry using dynamic programming (\cite{Riddler_Question} \cite{Riddler_Answer} \cite{GitHub}).  Our hope is to use the graphs described in this paper to develop a combinatorial enumeration of crossword grids which does not rely on computer assistance.

\bibliographystyle{plain}
\bibliography{main}

\end{document}